\documentstyle[11pt]{article}

\begin{document}
\noindent
\begin{center}
  {\LARGE Orbifold Quantum Cohomology}\footnote{both authors partially
supported by the National Science Foundation}
  \end{center}

  \noindent
  \begin{center}

     {\large Weimin Chen\footnote{New Address (starting Fall 2000):
Department of Mathematics, SUNY at Stony Brook, NY 11794}
             and Yongbin Ruan}\\[5pt]
      Department of Mathematics, University of Wisconsin-Madison\\
        Madison, WI 53706\\[5pt]
        Email: wechen@math.wisc.edu  and ruan@math.wisc.edu\\[5pt]
              \end{center}

              \def \x{{\bf x}}
              \def \J{{\cal J}}
              \def \M{{\cal M}}
              \def \A{{\cal A}}
              \def \B{{\cal B}}
              \def \C{{\bf C}}
              \def \Z{{\bf Z}}
              \def \R{{\bf R}}
              \def \P{{\bf P}}
              \def \I{{\bf I}}
              \def \N{{\cal N}}
              \def \T{{\cal T}}
              \def \Q{{\bf Q}}
              \def \D{{\cal D}}
              \def \H{{\cal H}}
              \def \S{{\cal S}}
              \def \e{{\bf E}}
              \def \CP{{\bf CP}}
              \def \U{{\cal U}}
              \def \E{{\cal E}}
              \def \F{{\cal F}}
              \def \L{{\cal L}}
              \def \K{{\cal K}}
              \def \G{{\cal G}}
              \def \z{{\bf z}}
              \def \m{{\bf m}}
              \def \n{{\bf n}}
              \def \V{{\cal V}}
              \def \W{{\cal W}}

\vspace{1cm}

This is a research announcement on a theory of Gromov-Witten
invariants and quantum cohomology of symplectic or projective
orbifolds. Our project started in the summer of 98 where our
original motivation was to study the quantum cohomology under
singular flops in complex dimension three. In this setting, we
allow our three-fold to have terminal singularities which can be
deformed into a symplectic orbifold. We spent the second half of
98 and most of spring of 99 to develop the foundation of
Gromov-Witten invariants over orbifolds, including the key conceptual
ingredient --- the notion of good map. In the April of 99, we were
lucky to meet R. Dijkgraaf who explained to us the orbifold string theory and
the role of twisted sectors. The twisted sector provides the
precise topological framework for our orbifold quantum
cohomology. Our theory of orbifold quantum cohomology was
virtually completed in the summer of 99. Here, we give an overlook
of the foundation of orbifold quantum cohomology while we leave its
applications in other fields such as birational geometry for future
research. The first part of our work has already appeared in
\cite{CR1}. The second part \cite{CR2} will appear shortly. We would like to
thank R. Dijkgraaf for bringing orbifold string theory to our
attention at the critical moment of our project. We are also
benefited from many discussions with  E. Zaslow. Finally, the
second author would like to thank E. Witten for many stimulating
discussions about orbifold string theory.

\section{Orbifold String Theory}

In 1985, Dixon, Harvey, Vafa and Witten considered string
compactification on Calabi-Yau orbifolds (arising as global
quotients $X/G$ by a finite group $G$) for the purpose of symmetry
breaking \cite{DHVW}. Later on, orbifold string theory became an
important part of string theory. For example, orbifolds provide
some of the simplest nontrivial models in string theory. Until
very recently, only physical argument for mirror symmetry had been
given for orbifold models \cite{GP}. In fact, orbifolds are such a
popular topic in string theory that a search on hep-th yields more
than 200 papers whose title contains orbifold. The reason that
orbifold string theory is interesting mathematically is that it
contains information which we do not have in the smooth case.
Roughly speaking, to have a consistent string theory, string
Hilbert space has to contain factors called twisted sectors.
Twisted sectors can be viewed as the contribution from
singularities. All other quantities such as correlation functions
have to contain the contribution from the twisted sectors. So far,
the twisted sectors are best understood in the context of
conformal field theory. It is our intention to initiate a program
to investigate the new geometry and topology of orbifolds caused
by the inclusion of twisted sectors. We should mention that the
orbifold string theory construction has only been carried out for
global quotients. However, it is well-known that most of Calabi-Yau
orbifolds are not global quotients. It seems to  also be important
to be able to construct orbifold string theory over general
Calabi-Yau orbifolds. Our orbifold quantum cohomology theory works
over arbitrary symplectic or projective orbifolds. We hope that our
construction
will shed some light on the construction of orbifold string theory
for general Calabi-Yau orbifolds.

   An orbifold, by definition, is a singular space. One can try to
desingularize a Calabi-Yau orbifold by the means of resolution or
deformation. To preserve the Calabi-Yau condition, we have to
restrict ourselves to the so-called crepant resolutions. It is natural
to ask for the relation between orbifold string theory and string
theory of its desingularization. In fact, this link provides some
of the most interesting mathematics from orbifold string theory.
In physics, orbifold string theory and ordinary string theory of its crepant
resolution appear to be two members in a family of theories. This strongly
suggests that there must be a relation between them. The
strongest predication is that they are the same. Indeed, this is
what physicists hope for. For quantum cohomology, this translates
into the following orbifold string theory prediction:

\vspace{2mm}

\noindent{(1.1)}\hspace{8mm}{\it The quantum cohomology of a
crepant resolution should be ``isomorphic'' to the `` quantum
cohomology'' of the orbifold.}

\vspace{2mm}

Here, the ``quantum cohomology'' of the orbifold should be understood
as orbifold quantum cohomology. The goal of our project is to
establish a mathematical theory for orbifold quantum cohomology.

    Actually, the above prediction is false in general. A counterexample is
a $K3$-surface with ADE-singularities. But this is not the end of
the  story. Recall that the most general form of mirror symmetry
fails for rigid Calabi-Yau 3-folds. But this did not stop research
from unearthing  the layer and layer of mathematical treasures
from mirror symmetry. In fact, it is entirely possible that weaker
forms or the current form $(1.1)$ for a more restrictive class of
orbifolds are still true. The authors believe that this link to
crepant resolutions will greatly enrich  this subject. Therefore,
it is useful to keep this strongest form of prediction in mind for
the direction of future research.

    The weakest form of the orbifold string theory prediction is to replace
quantum cohomology by orbifold Euler number. Here, it has a good
chance to hold. The orbifold Euler number is defined as the sum of
Euler numbers over all sectors. It has a natural interpretation as
the Euler characteristic of orbifold K-theory \cite{AR}. A weak
orbifold string theory prediction is that Euler number of a
crepant resolution is the same as the orbifold Euler number of
itself. On the mathematical side, a similar phenomenon was
independently discovered earlier by John McKay, which is now known
as {\it McKay correspondence} \cite {Mc,Re}. A version of
McKay correspondence is stated as follows:

{\it Let $G\subset SL(n,\C)$ be a finite group, and
$\pi:Y\rightarrow X=\C^n/G$ be a crepant resolution, then there
exist ``natural'' bijections between conjugacy classes of $G$ and
a basis of $H_{\ast}(Y;\Z)$.}

Based on these ideas, Batyrev-Dais \cite{BD} proposed the
so-called {\it strong McKay correspondence} and defined {\it
string-theoretic Hodge numbers}.

The classical part of our orbifold quantum cohomology is a new cohomology
of orbifolds which we call orbifold cohomology (see section 2). In the
case of Gorenstein orbifolds, Batyrev-Dais's string-theoretic
Hodge number is just the Hodge number of our orbifold cohomology.
The next level of the orbifold string theory prediction is to
identify the orbifold cohomology group with the ordinary
cohomology group of a crepant resolution. This is best described
through the orbifold K-theory \cite{AR}. It is unlikely that one
can identify the cohomology ring structures because of quantum
corrections. The third level would be the last level concerning
quantum cohomology, which is the most challenging one. At this moment,
it is not clear how to
formulate the prediction without the risk of finding a simple
counterexample.

\vspace{2mm}

In the following sections, we shall outline the construction of an
{\it orbifold Gromov-Witten theory}, which obeys almost all of the
axioms of ordinary Gromov-Witten theory, as it should according to
physics.

\section{Orbifold Cohomology}

 The ordinary quantum
cohomology ring appears as a (quantum) deformation of the ordinary
cohomology ring with the cup product. In the orbifold
Gromov-Witten theory, the role of ordinary cohomologies is played
by the so called {\it orbifold cohomologies}, which we shall
describe in this section. One of our main results is the
construction of an {\it orbifold cup product} on the total orbifold
cohomology group, which makes it into a ring with unit. We will
call the resulting ring {\it the orbifold cohomology ring}. The
orbifold quantum cohomology ring is just the corresponding quantum
deformation of the orbifold cohomology ring. Details are in
\cite{CR1}.

\vspace{2mm}

Let $X$ be a closed, almost complex orbifold of dimension $n$,
with almost complex structure $J$. For any $p\in X$, let
$(V_p,G_p,\pi_p)$ be a uniformizing system, which can be taken so
that $V_p$ is a small ball in $\C^n$ centered at the origin and
$G_p$ acts on $V_p$ as a finite subgroup of $U(n)$. We consider
the set $$ \widetilde{X}:=\{(p,(g)_{G_p})|p\in X, g\in G_p\}.
\leqno(2.1) $$ Here $(g)_{G_p}$ represents the conjugacy class of
$g$ in $G_p$. There is a locally constant function
$\iota:\widetilde{X}\rightarrow\Q$ defined as follows: write $g$
as a diagonal matrix $$ diag(e^{2\pi i m_{1,g}/m_g}, \cdots,
e^{2\pi i m_{n,g}/m_g}), $$ where $m_g$ is the order of $g$ in
$G_p$, and $0\leq m_{i,g}<m_g$ for $i=1,\cdots,n$. We define $$
\iota(p,(g)_{G_p})=\sum_{i=1}^n\frac{m_{i,g}}{m_g}.\leqno(2.2) $$
Let $I:\widetilde{X}\rightarrow\widetilde{X}$ be the involution
defined by $I(p,(g)_{G_p})=(p,(g^{-1})_{G_p})$.

\vspace{3mm}

\noindent{\bf Lemma 2.1:}\hspace{2mm}{\it There is an equivalence
relation among the $(g)_{G_p}$, and if we let $T=\{(g)\}$ be the
set of equivalence classes and define $X_{(g)}=\{(p,(g)_{G_p})\in
\widetilde{X}|(g)_{G_p}\in (g)\}$, then each $X_{(g)}$ is
naturally a closed, connected, almost complex orbifold, and
$\widetilde{X}$ is decomposed as a disjoint union $\sqcup_{(g)\in
T}X_{(g)}$. Furthermore, if we denote the value of the locally
constant function $\iota:\widetilde{X}\rightarrow\Q$ by
$\iota_{(g)}$, and let $(g^{-1})$ denote the image of $(g)$ under
the involution $I$, and $(1)$ denote the equivalence class of the
trivial element $(1)_{G_p}$, the following conditions are
satisfied:
\begin{enumerate}
\item $\iota:\widetilde{X}\rightarrow\Q$ is integer-valued iff each $G_p$
is contained in $SL(n,\C)$.
\item $$
\iota_{(g)}+\iota_{(g^{-1})}=\dim_\C X-\dim_\C X_{(g)}.\leqno(2.3)
$$
\item $\iota_{(g)}\geq 0$ for all $(g)\in T$, and $\iota_{(g)}=0$ iff
$(g)=(1)$.
\end{enumerate}
}
Note that for Calabi-Yau orbifolds, each $\iota_{(g)}$ is integer-valued.
When $X=Y/G$ is a global quotient, $\widetilde{X}$ can be identified with
$\sqcup_{(g),g\in G}Y^g/C(g)$.

The orbifold cohomologies are just direct sums of ordinary cohomologies of
$X_{(g)}$ with degrees shifted by $2\iota_{(g)}$. More precisely,

\vspace{2mm}

\noindent{\bf Definition 2.2:}\hspace{3mm} {\it Let $X$ be a
closed almost complex orbifold with $\dim_{\C} X=n$. For any
rational number $d\in [0,2n]$, the orbifold cohomology group of
degree $d$ is defined to be the direct sum $$
H^d_{orb}(X;\Q)=\oplus_{(g)\in
T}H^{d-2\iota_{(g)}}(X_{(g)};\Q).\leqno(2.4) $$ }
We will call
$\iota_{(g)}$ {\it degree shifting numbers}, which have been
referred as {\it fermion shift numbers} in physics \cite{Z}. The
orbifold $X_{(g)}$ or its cohomology will be called a {\it twisted
sector} if $(g)\neq (1)$, and called {\it the nontwisted sector}
if $(g)=(1)$. The construction of $\widetilde{X}$ (cf. (2.1))
first appeared in \cite{Ka}.

\vspace{2mm}

The following theorem is proved in \cite{CR1}, whose construction
is based on genus-zero, degree zero orbifold Gromov-Witten
invariants.

\vspace{2mm}

\noindent{\bf Theorem 2.3:}\hspace{2mm}{\it
Let $(X,J)$ be a closed almost complex orbifold of dimension $n$. Then
\begin{enumerate}
\item There is a non-degenerate pairing $<,>_{orb}:
H^d_{orb}(X;\Q)\times H^{2n-d}_{orb}(X;\Q)\rightarrow \Q$ extending the
ordinary Poincar\'{e} pairing on the nontwisted sectors $H^\ast(X;\Q)$.
\item There is a cup product $\cup_{orb}: H_{orb}^p(X;\Q)\times
H_{orb}^q(X;\Q) \rightarrow H_{orb}^{p+q}(X;\Q)$ for any $0\leq
p,q\leq 2n$ such that $p+q\leq 2n$, which has the following
properties:
\begin{itemize}
\item The total orbifold cohomology group $H^\ast_{orb}(X;\Q)=
\oplus_{0\leq d\leq 2n}H^d_{orb}(X;\Q)$ is a ring with unit
$e_X^0\in H_{orb}^0(X;\Q)$ under $\cup_{orb}$, where $e_X^0$ is
the Poincar\'{e} dual to the fundamental class $[X]$.
In particular, $\cup_{orb}$ is associative.
\item Restricted to each $H_{orb}^d(X;\Q)\times H_{orb}^{2n-d}(X;\Q)
\rightarrow H_{orb}^{2n}(X;\Q)$,
$$
\int_X \alpha\cup_{orb}\beta=<\alpha, \beta>_{orb}. \leqno(2.5)
$$
\item The cup product $\cup_{orb}$ is invariant under deformations of $J$.
\item When $X$ is of integral degree shifting numbers, the total orbifold
cohomology group $H_{orb}^\ast(X;\Q)$ is integrally graded, and we
have supercommutativity
$$
\alpha_1\cup_{orb}\alpha_2=(-1)^{\deg\alpha_1\cdot\deg\alpha_2}\alpha_2
\cup_{orb}\alpha_1. \leqno(2.6)
$$
\item Restricted to the nontwisted sectors, i.e., the ordinary
cohomologies $H^\ast(X;\Q)$, the cup product $\cup_{orb}$ equals the
ordinary cup product on $X$.
\end{itemize}
\end{enumerate}
}

We remark that there is an analogous construction using Dolbeault
cohomology groups; for details see \cite{CR1}.

\section{Good Map and Pull-Back Bundle}

Now we come to one of the main issues in the construction of
orbifold Gromov-Witten invariants. Recall that if $f:X\rightarrow
X^\prime$ is a $C^\infty$ map between manifolds and $E$ is a
smooth vector bundle over $X^\prime$, then there is a smooth
pull-back vector bundle $f^\ast E$ over $X$ and a bundle morphism
$\bar{f}:f^\ast E\rightarrow E$ which covers the map $f$. However,
if instead we have a $C^\infty$ map $\tilde{f}$ between orbifolds
$X$ and $X^\prime$, and an orbibundle $E$ over the orbifold
$X^\prime$, the question whether there is a pull-back orbibundle
$E^\ast$ over $X$ and an orbibundle morphism
$\bar{f}:E^\ast\rightarrow E$ covering the map $\tilde{f}$ is a
quite complicated issue: $E^\ast$ might not exist, or even if it
exists it might not be unique. Traditionally, a neighborhood of a
smooth map into a manifold is described by smooth sections of the
pull-back of the tangent bundle of the manifold. Hence
understanding this question became the very first step in
describing the moduli spaces of pseudo-holomorphic maps, or more
precisely, the very first step in order to understand what would
be the corresponding notion of ``stable map'' in the orbifold
case.

By a {\it $C^\infty$ map} between orbifolds $X$ and $X^\prime$ we
mean an equivalence class of collections of local smooth liftings
between local uniformizing systems of a continuous map from $X$ to
$X^\prime$. This notion is equivalent to the notion of {\it
V-maps} in \cite{S}, where the notion of orbifold was first
introduced under the name {\it V-manifold}. A brief review of
orbifolds is given in \cite{CR1}, and a self-contained, elementary
discussion of various aspects of differential geometry and global
analysis on orbifolds is contained in \cite{CR2}.

\vspace{2mm}

Now we describe our key concept: the notion of {\it good map}.
Let $\tilde{f}:X\rightarrow X^\prime$ be a $C^\infty$ map between orbifolds $X$
and $X^\prime$ whose underlying continuous map is denoted by $f$.
Let $\U=\{U_\alpha;\alpha\in\Lambda\}$ be an open cover of $X$ and
$\U^\prime=\{U^\prime_\alpha;\alpha\in\Lambda\}$ be an open cover of the
image $f(X)$ in $X^\prime$, which satisfy the following conditions:

\begin{enumerate}
\item Each $U_\alpha$ (resp. $U^\prime_\alpha$) is uniformized by
$(V_\alpha,G_\alpha,\pi_\alpha)$ (resp. by $(V^\prime_\alpha,
G^\prime_\alpha,\pi^\prime_\alpha)$).
\item If $U_\alpha\subset U_\beta$ (resp.
$U^\prime_\alpha\subset U^\prime_\beta$), then there is a collection of
smooth open embeddings $(V_\alpha,G_\alpha,\pi_\alpha)\rightarrow
(V_\beta,G_\beta,\pi_\beta)$ (resp. $(V^\prime_\alpha,G^\prime_\alpha,
\pi^\prime_\alpha)\rightarrow (V^\prime_\beta,G^\prime_\beta,
\pi^\prime_\beta)$), which are called {\it injections}.
\item For any point $p\in U_\alpha\cap U_\beta$ (resp. $p^\prime
\in U^\prime_\alpha\cap U^\prime_\beta$), there is a $U_\gamma$ (resp. $U^\prime_\gamma$) such that $p\in U_\gamma\subset U_\alpha\cap U_\beta$ (resp.
$p^\prime\in U^\prime_\gamma\subset U^\prime_\alpha\cap U^\prime_\beta$).
\item Any inclusion $U_\alpha\subset U_\beta$ implies $U^\prime_\alpha
\subset U^\prime_\beta$.
\item For each $\alpha\in\Lambda$, $f(U_\alpha)\subset U^\prime_\alpha$.
Moreover, there is a collection of local smooth liftings of $f$,
$\{\tilde{f}_\alpha: V_\alpha\rightarrow V^\prime_\alpha;\alpha\in\Lambda\}$
which defines the given $C^\infty$ map $\tilde{f}$, such that any
injection $i_{\beta\alpha}:(V_\alpha,G_\alpha,\pi_\alpha)\rightarrow
(V_\beta,G_\beta,\pi_\beta)$ is assigned with an injection
$\lambda(i_{\beta\alpha}):(V^\prime_\alpha,G^\prime_\alpha,
\pi^\prime_\alpha)\rightarrow (V^\prime_\beta,G^\prime_\beta,
\pi^\prime_\beta)$ satisfying the following compatibility conditions:
$$
\tilde{f}_\beta\circ i_{\beta\alpha}
=\lambda(i_{\beta\alpha})\circ\tilde{f}_\alpha
\hspace{3mm}\forall \alpha,\beta\in\Lambda
\leqno(3.1)
$$ and
$$
\lambda(i_{\gamma\beta}\circ i_{\beta\alpha})=\lambda(i_{\gamma\beta})\circ
\lambda(i_{\beta\alpha}) \hspace{3mm}\forall \alpha,\beta,\gamma\in\Lambda.
\leqno(3.2)
$$
\end{enumerate}

\noindent{\bf Definition 3.1:}\hspace{2mm}{\it We call such a
$(\U,\U^\prime,\{\tilde{f}_\alpha\},\lambda)$ a compatible system of
the $C^\infty$ map $\tilde{f}$.
A $C^\infty$ map is said to be good if it admits a compatible system.}

\vspace{2mm}

\noindent{\bf Lemma 3.2:}\hspace{2mm}{\it
Let $pr:E\rightarrow X^\prime$ be an orbibundle over $X^\prime$.
For any $C^\infty$ good map $\tilde{f}:X\rightarrow X^\prime$ with a
compatible system $\xi=(\U,\U^\prime,\{\tilde{f}_\alpha\},\lambda)$,
there is a canonically defined pull-back orbibundle
$pr:E^\ast_{\tilde{f},\xi}\rightarrow X$
with an orbibundle morphism $\bar{f}_{\tilde{f},\xi}:
E^\ast_{\tilde{f},\xi}\rightarrow E$
which covers $\tilde{f}$.
}

\vspace{2mm}

\noindent{\bf Definition 3.3:}\hspace{2mm}{\it Two compatible
systems $\xi_1,\xi_2$ of a $C^\infty$ good map
$\tilde{f}:X\rightarrow X^\prime$ are said to be isomorphic if for
any orbibundle $E$ over $X^\prime$ there is an orbibundle
isomorphism $\phi:E^\ast_{\tilde{f},\xi_1}\rightarrow
E^\ast_{\tilde{f},\xi_2}$ such that $$
\bar{f}_{\tilde{f},\xi_1}=\bar{f}_{\tilde{f},\xi_2}\circ\phi.\leqno(3.3)
$$ }

\vspace{2mm}

\noindent {\bf Example 3.4a: }{\it Not every $C^\infty$ map is
good, as shown in the following example: consider an effective
linear representation of a finite group $(\R^n,G)$. Let $H^g$ be
the linear subspace of fixed points of an element $1\neq g\in G$.
Then the centralizer $C(g)$ of $g$ in $G$ acts on $H^g$, and
$(H^g,C(g)/K_g)$ is an effective linear representation, where
$K_g\subset C(g)$ is the kernel of the action of $C(g)$ on $H^g$.
Suppose $H^g\neq \{0\}$ and there is no homomorphism $\lambda:
C(g)/K_g\rightarrow C(g)$ such that $\pi\circ\lambda$ is the
identity homomorphism, where $\pi: C(g)\rightarrow C(g)/K_g$ is
the projection, then the continuous map $H^g/C(g)\rightarrow \R^n/G$ induced
by inclusion $H^g\hookrightarrow \R^n$ is a
$C^\infty$ map which is not a good one. }

\vspace{2mm}

\noindent{\bf Example 3.4b: }{\it
There could be non-isomorphic compatible systems of
the same $C^\infty$ map, as shown in the following example:
Let $X=\C\times\C/G$ where $G=\Z_2
\oplus\Z_2$ acting on $\C\times\C$ in the standard way. For the $C^\infty$
map $\C/\Z_2\rightarrow X$ defined by the
inclusion $\C\times\{0\}\hookrightarrow \C\times\C$,
there are two non-isomorphic compatible systems $(\tilde{f}_i,\lambda_i):
(\C,\Z_2)\rightarrow (\C\times\C,G)$, $i=1,2$, where $\lambda_1(1)=(1,0)$
and $\lambda_2(1)=(1,1)$.
}

\vspace{2mm}

In general it is not only difficult to determine whether a given
$C^\infty$ map is good or not, but also difficult to classify
compatible systems of an arbitrary good map by equivalence up to
isomorphism. Nevertheless, it is clear that a good map together
with an isomorphism class of compatible systems is the object we
ought to deal with in the orbifold quantum cohomology theory.

\vspace{2mm}

We end this section with a discussion on the case when the domain of a good
map is a 2-dimensional orbifold.

\vskip 0.1in

\noindent{\bf Definition-Construction 3.5:}\hspace{2mm}{\it
\begin{enumerate}
\item Given a marked Riemann surface $(\Sigma,\z)$ where $\z=(z_1,\cdots,z_k)$
is the set of marked points, we can give a unique orbifold
structure to $\Sigma$ by assigning to each marked point $z_i$  an
integer $m_i\geq 1$ (note that $m_i=1$ is allowed for
convenience). We will call $(\Sigma,\z,\m)$ an orbifold marked
Riemann surface, where $\m=(m_1,\cdots,m_k)$ is the set of
assigned integers, called multiplicities.
\item An orbifold nodal Riemann surface is a marked nodal Riemann surface with
the following data: (i) each irreducible component is an orbifold marked
Riemann surface (here a nodal point is considered marked on an irreducible
component); (ii) two identified nodal points are assigned with the same
multiplicity.
\end{enumerate}
}

\vspace{2mm}

\noindent{\bf Convention-Definition 3.6:}\hspace{2mm}{\it
\begin{enumerate}
\item Note that, in the definition of compatible systems, the compatibility
conditions $(3.1),(3.2)$ give rise to a collection of
homomorphisms $\lambda_\alpha, \alpha\in\Lambda$, between local
groups $G_\alpha$ and $G^\prime_\alpha$, such that each local
smooth lifting $\tilde{f}_\alpha:V_\alpha\rightarrow
V_\alpha^\prime$ is $\lambda_\alpha$-equivariant. For a good
$C^\infty$ map whose domain is an orbifold marked Riemann surface,
we require that each $\lambda_\alpha$ be a monomorphism for any of
its compatible systems.
\item A good $C^\infty$ map with a compatible system from an orbifold nodal
Riemann surface into an orbifold $X$ is a collection of good
$C^\infty$ maps with compatible systems defined on its irreducible
components which satisfies the following compatibility condition:
for each pair of identified nodal points $z_\nu$ and $z_\omega$,
the homomorphisms $\lambda_{z_\nu}$ and $\lambda_{z_\omega}$
between local groups, which are determined by the corresponding
compatible systems, satisfy the equation $$
\lambda_{z_\nu}(x)\cdot\lambda_{z_\omega}(x)=1_{G_p} \leqno(3.4)
$$ in $G_p$, where $p\in X$ is the image of the identified nodal
points $z_\nu$ and $z_\omega$ under the good $C^\infty$ map, and
$x$ is a generator of the local cyclic group at $z_\nu$ and
$z_\omega$ ($z_\nu$ and $z_\omega$ have the same multiplicity, hence
the same local cyclic group).
\end{enumerate}
}

\vskip 0.1in

Finally we observe that each good $C^\infty$ map with a compatible
system from an orbifold nodal Riemann surface with $k$ marked
points into an orbifold $X$ determines a point in the space
$\widetilde{X}^k$,  where $\widetilde{X} =\sqcup_{(g)\in
T}X_{(g)}$ (cf. (2.1)), as follows: let the underlying continuous
map be $f$ and for each marked point $z_i$, $i=1,\cdots,k$, let
$x_i$ be the positive generator of the cyclic local group at
$z_i$, and $\lambda_{z_i}$ be the homomorphism determined by the
given compatible system, then the determined point in $\widetilde{X}^k$ is
$$
((f(z_1),(\lambda_{z_1}(x_1))_{G_{f(z_1)}}),\cdots,
(f(z_k),(\lambda_{z_k}(x_k))_{G_{f(z_k)}})).\leqno(3.5)
$$
Let
$\x=(X_{(g_1)},\cdots,X_{(g_k)})$ be a connected component in
$\widetilde{X}^k$. We say that a good map with a compatible system
is of {\it type} $\x$ if the point $(3.5)$ it determines in
$\widetilde{X}^k$ lies in the component $\x$.

\section{Orbifold Stable Maps}

We start with the definition of {\it pseudo-holomorphic map} from a Riemann
surface into an almost complex orbifold.

\vspace{2mm}

\noindent{\bf Definition 4.1:}\hspace{2mm}{\it
A pseudo-holomorphic map from a Riemann surface $(\Sigma,j)$ into an almost
complex orbifold $(X,J)$ is a continuous map $f: \Sigma\rightarrow X$ which
satisfies the following conditions:
\begin{enumerate}
\item For any point $z\in \Sigma$, there is a disc neighborhood $D_z$ of $z$
with a branched covering map $br_z:\widetilde{D}_z\rightarrow D_z$ given by
$w\rightarrow w^{m_z}$ (here $m_z=1$ is allowed).
\item Let $p=f(z)$. There is a local uniformizing system
$(V_p,G_p,\pi_p)$ of $X$ at $p$ and a local smooth lifting $\tilde{f}_z:
\widetilde{D}_z\rightarrow V_p$ of $f$ in the sense that
$f\circ br_z=\pi_p\circ\tilde{f}_z$.
\item $\tilde{f}_z$ is pseudo-holomorphic, i.e.,
$d\tilde{f}_z\circ j=J\circ d\tilde{f}_z$.
\end{enumerate}
}

\vspace{2mm}

\noindent{\bf Remarks 4.2:}{\it
\begin{enumerate}
\item When $(X,J)$ is a complex orbifold, i.e., $J$ is integrable, a
pseudo-holomorphic map $f:(\Sigma,j)\rightarrow (X,J)$ is just a holomorphic
map from $(\Sigma,j)$ into the analytic space $(X,J)$.
\item For each pseudo-holomorphic map $f:(\Sigma,j)\rightarrow (X,J)$, there
is a subset of finitely many points
$\{z_1,z_2,\cdots,z_k\}\subset\Sigma$ such that for any
$z\in\Sigma\setminus\{z_1,z_2,\cdots,z_k\}$ the multiplicity $m_z$
in Definition 4.1-1 equals one (cf. \cite{HW}). We will consider
pseudo-holomorphic maps from marked Riemann surfaces into $(X,J)$.
As a convention we will always mark these points
$\{z_1,z_2,\cdots,z_k\}$ where the multiplicity is greater than
one.
\end{enumerate}
}

\vspace{2mm}

Given a pseudo-holomorphic map $f$ from a marked Riemann surface
$(\Sigma,\z)$ into $(X,J)$, where $\z=(z_1,\cdots,z_k)$ is a set
of finitely many distinct marked points on $\Sigma$, there is an
orbifold structure on $\Sigma$ with singular set contained in $\z$
such that $f$ can be lifted to a good $C^\infty$ map $\tilde{f}$.
A crucial technical result is summarized in the following

\vspace{2mm}

\noindent{\bf Lemma 4.3:}\hspace{2mm}{\it For any
pseudo-holomorphic map $f$ from a Riemann surface $\Sigma$ of
genus $g$ with $k$ marked points $\z=(z_1,z_2,\cdots,z_k)$ into
$(X,J)$, there are finitely many orbifold structures on $\Sigma$
whose singular set is contained in $\z$, and for each of these
orbifold structures there are finitely many pairs
$(\tilde{f},\xi)$, where $\tilde{f}$ is a good map whose
underlying map is $f$, and $\xi$ is an isomorphism class of
compatible systems of $\tilde{f}$. The total number is bounded
from above by a constant $C(X,g,k)$ depending only on $X,g,k$. }

\vskip 0.1in

\noindent{\bf Definition 4.4:}\hspace{2mm}{\it
An orbifold stable map from a marked nodal Riemann surface into an almost
complex orbifold $(X,J)$ consists of the following data:
\begin{enumerate}
\item A continuous map from the marked nodal Riemann surface into $(X,J)$
whose restriction to each irreducible component is pseudo-holomorphic.
\item An orbifold structure on the marked nodal Riemann surface so that it
becomes an orbifold nodal Riemann surface, and a good map with a compatible
system from the orbifold nodal Riemann surface into $(X,J)$ with the given
underlying continuous map.
\item Stability condition: on each $S^2$ or $T^2$ component which is mapped into
a point in $X$ there are at least three or one special points
(marked or nodal).
\end{enumerate}
There is an obvious equivalence
relation amongst the set of orbifold stable maps. We denote by
$\overline{\M}_{g,k}(X,J,A,\x)$ the set of all equivalence classes of
orbifold stable maps of
genus $g$, $k$ marked points, type $\x$ and homology class $A$ into $(X,J)$.
}

\vspace{2mm}

\noindent {\bf Remark 4.5: }{\it In the algebraic setting of
Deligne-Mumford stack, a related notion which is called twisted
stable map was discussed in \cite{AV}. Their twisted stable map
was described in the language of category and functor. Our good
map was formulated in elementary differential-geometric language.
From the first sight, two notions look quite different. However,
D. Abramovich kindly informed us that they are actually equivalent
\cite{A}.
 }

\vskip 0.1in

\noindent {\bf Remark 4.6: }{\it If $f: \Sigma\rightarrow X$ is a
pseudo-holomorphic map whose image intersects the singular locus
of $X$ at only finitely many points, then there is a unique choice
of orbifold structure on $\Sigma$ together with a unique
$(\tilde{f},\xi)$, where $\tilde{f}$ is a good map with an
isomorphism class of compatible systems $\xi$ whose underlying
continuous map is $f$. If the image of $f$ lies completely inside
the singular locus, there could be different choices, and they are
regarded as different points in the moduli space. }

\vspace{2mm}

\noindent{\bf Definition 4.7:}\hspace{2mm}{\it
\begin{enumerate}
\item An orbifold $X$ is symplectic if there is a closed 2-form $\omega$ on
$X$ whose local liftings are non-degenerate.
\item A projective orbifold is a complex orbifold which is a projective
variety as an analytic space.
\end{enumerate}
}

\vspace{2mm}

The usual Gromov Compactness Theorem for pseudo-holomorphic maps combined with
Lemma 4.3 gives the following

\vspace{2mm}

\noindent{\bf Proposition 4.8:}\hspace{2mm}{\it Suppose that $X$
is a symplectic or projective orbifold. The moduli space of
orbifold stable maps $\overline{\M}_{g,k}(X,J,A, \x)$ is a compact
metrizable space under a natural topology, whose ``virtual
dimension'' is $2d$, where $$ d=c_1(TX)\cdot
A+(\dim_{\C}X-3)(1-g)+k-\iota(\x).\leqno $$ Here
$\iota(\x):=\sum_{i=1}^k \iota_{(g_i)}$ for
$\x=(X_{(g_1)},\cdots,X_{(g_k)})$. }

\vspace{2mm}

\section{Orbifold Quantum Cohomology}

For any component $\x=(X_{(g_1)},\cdots,X_{(g_k)})$, there are $k$
evaluation maps (cf. (3.5)) $$
e_i:\overline{\M}_{g,k}(X,J,A,\x)\rightarrow X_{(g_i)},
\hspace{4mm} i=1,\cdots, k. \leqno(5.1) $$
For any set of cohomology classes
$\alpha_i\in H^{*-2\iota_{(g_i)}}(X_{(g_i)};\Q)\subset
H^*_{orb}(X;\Q)$, $i=1,\cdots,k$, the orbifold Gromov-Witten
invariant is defined as the
virtual integral $$ \Psi^{X,J}_{(g,k,A,\x)}(\alpha^{l_1}_1, \cdots,
\alpha^{l_k}_k)=\int^{vir}_{
\overline{\M}_{g,k}(X,J,A,\x)}\prod_{i=1}^k c_1(L_i)^{l_i}e^*_i
\alpha_i,\leqno(5.2) $$ where $L_i$ is the line bundle generated
by cotangent space of the $i$-th marked point.

\vspace{2mm}

When $g=0$ and $A=0$, the moduli space
$\overline{\M}_{g,k}(X,J,A,\x)$ admits a very nice and elementary
description, based on which we gave an elementary construction of
genus zero, degree zero orbifold Gromov-Witten invariants in \cite{CR1}.
Even in this case, virtual integration is needed where there is an
obstruction bundle. The orbifold cup product (cf. Theorem 2.3) is
defined through these orbifold Gromov-Witten invariants. In the
general case, we need to use the full scope of the virtual
integration machinary developed by \cite{FO}, \cite{LT}, \cite{Ru}
and \cite{Sie}.

\vspace{2mm}

Singularities of an orbifold impose additional difficulties in
carrying  out virtual integration  in the orbifold case. Due to
the presence of singularities, even on a closed orbifold, the
function of injective radius of the exponential map does not have
a positive lower bound. As a consequence, it is not known that a
neighborhood of a (good) $C^\infty$ map into an orbifold can be
completely described by $C^\infty$ sections of the pull-back
tangent bundle via the exponential map. Our approach is a
combination of techniques developed in the smooth case. We first
construct a local Kuranishi neighborhood for each stable map in
$\overline{\M}_{g,k}(X,J,A,\x)$, then find finitely many stable
maps whose local Kuranishi neighborhoods (although they may have
different dimensions) can be patched together to form a ``global
virtual neighborhood'' of $\overline{\M}_{g,k}(X,J,A,\x)$ (cf.
\cite{FO}). This is similar to the constructions of \cite{FO},
\cite{LT}. We carry out the virtual integration over this ``global
virtual neighborhood'' by constructing a system of compatible
``Thom forms'' (cf. \cite{Ru}).

When $X$ has a symplectic torus action, the ``global virtual
neighborhood'' can be constructed so that it respects this torus
action. The localization theory can be extended to the case of
virtual integration. We leave this to another paper.

\vspace{2mm}

Main results of this work are summarized in the following

\vspace{2mm}

\noindent{\bf Theorem 5.1:}\hspace{2mm}{\it Let $X$ be a closed
symplectic or projective orbifold. The orbifold Gromov-Witten
invariants defined in (5.2) satisfy the quantum cohomology axioms
of Witten-Ruan for ordinary Gromov-Witten invariants (cf.
\cite{Ru1}) except that in the Divisor Axiom, the divisor class
is required to be in the nontwisted sector (i.e. in $H^2(X;\Q)$).
In the formulation of axioms, the ordinary cup product is
replaced by the orbifold cup product $\cup_{orb}$ (cf. Theorem 2.3). }

\vspace{2mm}

As a consequence, we have

\vspace{2mm}

\noindent{\bf Theorem 5.2:}\hspace{2mm}{\it Let $X$ be a closed
symplectic or projective orbifold. With suitable coefficient ring
$\cal{C}$, the small quantum product and the big quantum product
are well-defined on $H^\ast_{orb}(X;\Q)\otimes\cal{C}$, and have
properties similar to those of the ordinary quantum cohomology. }

\section{Closing Remarks}
    What we have accomplished so far is just a tip of iceberg! For example,
it is still a difficult problem to compute orbifold quantum
cohomology. This requires developing new machinery such as
localization and surgery techniques. There are two topics whose
natural home should be orbifold. They are  birational geometry and
mirror symmetry. For birational geometry, recent results in
algebraic geometry show that birational transformation can be
decomposed as a sequence of wall-crossings in GIT-quotients
\cite{AKMW}, \cite{HK}, \cite{W1}, \cite{W2}. The latter is
naturally in the orbifold category. For mirror symmetry, the
Calabi-Yau 3-folds in most of the known examples are crepant
resolutions of Calabi-Yau orbifolds. Therefore, it is more natural
to consider mirror symmetry for orbifolds. Moreover, the second
author believes that orbifold quantum cohomology is different from
the quantum cohomology of crepant resolutions. How to formulate
mirror symmetry in the categary of Calabi-Yau orbifolds seems to
be an extremely interesting problem. Suppose we can do all of
these, we are still working only in the so-called type II string
theory. There are orbifold versions for other types of string
theory (such as heterotic string theory) as well. The amount of
new mathematics we can unearth is unimaginable!

\end{document}